\documentclass[twoside,12pt,reqno]{amsart}

\usepackage{amsmath,amsthm,amssymb,amstext,amsfonts,amscd}
\usepackage{graphicx}
\usepackage{multirow}
\usepackage{url}
\newtheorem{definition}{Definition}[section]
 \newtheorem{theorem}{Theorem}[section]

\newtheorem{remark}{Remark}[section]

\numberwithin{equation}{section}

\begin{document}

\title[ Bell based Apostol type polynomials and its properties]
{ Bell based Apostol type polynomials and its properties}

\author[{\bf N. U. Khan and S. Husain}]{\bf Nabiullah Khan and Saddam Husain}

\address{Nabiullah Khan: Department of Applied
Mathematics, Faculty of Engineering and Technology,
      Aligarh Muslim University, Aligarh 202002, India}
 \email{nukhanmath@gmail.com}

\bigskip
\address{Saddam Husain: Department of Applied
    Mathematics, Faculty of Engineering and Technology,
    Aligarh Muslim University, Aligarh 202002, India}
 \email{saddamhusainamu26@gmail.com}



\begin{abstract}
 Many authors studied mix type polynomials and their properties. For their numerous uses in Statistics, Combinatorial analysis and Number theory. In the present article, we define a generating function of mix type Bell based Apostol type polynomials and numbers of order $\alpha$. Furthermore, we derive some elementary properties and identities of Bell based Apostol type polynomials of order $\alpha$ such as correlation formulas, Implicit summation formulas, deriative formulas and their special cases. Also, we derive some relation with Stirling numbers.
\end{abstract}

\maketitle

\section{\bf{Introduction}}

Recently, many authors studied various mix type polynomials and their applications. U. Duran et al. (see\cite{Araci-Acikgoz}) study the Bell Based Bernaulli polynomials and their applications. N.U. khan et al. (see \cite{khan-husain}) study Bell based Euler polynomials and their application. Motivated by these mention work, in this paper we study Bell based Apostol type polynomials and its properties such as correlation formula, derivative formula, implicit summation formula and their special cases. Moreover, we define some relation with Stirling numbers. 

Throughout paper, the symbols $\mathbb{Z}$ denote set of integers, $\mathbb{N}$ denotes set of natural numbers, $\mathbb{N}_{0}$ denoted set of non negative integers,  $\mathbb{R}$ denotes set of real numbers and  $\mathbb{C}$ denotes set of complex numbers respectively.
Recently, many authors (see \cite{Bell,Carlitz,Araci-Acikgoz,Kim-Kim,Srivastva-garg,Luo-ming-sri,luo-ming,luo-ming2} )  have studied the bivariate Bell polynomials, classical Bell polynomials, Bernoulli polynomials and Euler polynomials and Genocchi polynomials define as follows:

The Bell polynomials in two variable $x_{1}, x_{2}$ such as  ${\mathfrak{B}_n}(x_{1},x_{2})$ are known as bivariate Bell polynomials. The generating functions of bivariate Bell polynomials describe as follows (see \cite{Araci-Acikgoz}):
\begin{equation}\label{2.1}
	\sum_{n\geq0}{\mathfrak{B}_n}(x_{1},x_{2})\frac{{t}^n}{n!}= e^{x_{1}{t}}\,e^{x_{2}(e^{{t}}-1)}.
\end{equation}

When we take $x_{1}=0$, ${\mathfrak{B}_n}(0;x_{2})={\mathfrak{B}_n}(x_{2})$ are called classical Bell polynomials (or exponential polynomials), which is describe by following generating function (see \cite{Bell,Boas,Carlitz,Kim-Kim,kim-kim-kim-kwon}) which is defined as follows:
\begin{equation}\label{2.2}
	\sum_{n\geq0}{\mathfrak{B}_n}(x_{2})\frac{{t}^n}{n!}= e^{x_{2}(e^{{t}}-1)}.
\end{equation}

If we take $x_{2}=1$ in \eqref{2.2} i.e. ${\mathfrak{B}_n}(0;1)={\mathfrak{B}_n}(1)= {\mathfrak{B}_n}$ are called Bell numbers which is defined as follows (see \cite{Bell,Boas,Carlitz,Kim-Kim,kim-kim-kim-kwon}):
\begin{equation}\label{2.3}
	\sum_{n\geq0}{\mathfrak{B}_n}\frac{{t}^n}{n!}= e^{(e^{{t}}-1)}.
\end{equation}

The Euler polynomials ${\mathsf{E}}_{n}(x_{1})$ and Euler numbers ${\mathsf{E}}_{n}(0)$ are defined by the following generating function to be

\begin{equation}\label{e1.3}
	e^{x_{1}t}\left(\frac{2}{e^{t}+1}\right)=\sum_{n=0}^{\infty}\frac{t^n}{n!}~~{\mathsf{E}}_{n}(x_{1}), \,\,\,\,\,\,\, ( |t|<\pi; 1^{\alpha}:=1 ).
\end{equation}

If we take $x_{1} = 0$ then the Euler numbers $\mathsf{E}_{n}(0):= \mathsf{E}_{n}$ are defined by

\begin{equation}\label{e1.4}
	\left(\frac{2}{e^{t}+1}\right)=\sum_{n=0}^{\infty}\frac{t^n}{n!}~\mathsf{E}_{n}, \,\,\,\,\,\,\, ( |t|<\pi; 1^{\alpha}:=1 ).
\end{equation}

The Bernoulli polynomials and Bernoulli numbers is introduced by Datolli et al. (see\cite{Dattoli}), which is defined by the following generating function to be

\begin{equation}\label{e1.5}
	e^{x_{1}t}\left(\frac{t}{e^{t}-1}\right)=\sum_{n=0}^{\infty}\frac{t^n}{n!}{\mathsf{B}}_{n}(x_{1}), \,\,\,\,\,\,\, ( |t|<2\pi; 1^{\alpha}:=1).
\end{equation}

If we take $x_{1} = 0$ then the Euler numbers ${\mathsf{B}}_{n}(0):= {\mathsf{B}}_{n}$ are defined by

\begin{equation}\label{e1.6}
	\left(\frac{t}{e^{t}-1}\right)=\sum_{n=0}^{\infty}\frac{t^n}{n!}{\mathsf{B}}_{n}, \,\,\,\,\,\,\, ( |t|<2\pi; 1^{\alpha}:=1 ).
\end{equation}

The Genocchi polynomials ${\mathsf{G}}_{n}(x_{1})$ and Euler numbers ${\mathsf{G}}_{n}(0)$ are defined by the following generating function to be (see \cite{Srivastva-garg})

\begin{equation}\label{e1.3}
	e^{x_{1}t}\left(\frac{2t}{e^{t}+1}\right)=\sum_{n=0}^{\infty}\frac{t^n}{n!}~~{\mathsf{G}}_{n}(x_{1}), \,\,\,\,\,\,\, ( |t|<\pi; 1^{\alpha}:=1 ).
\end{equation}

If we take $x_{1} = 0$ then the Genocchi numbers $\mathsf{G}_{n}(0):= \mathsf{G}_{n}$ are defined by

\begin{equation}\label{e1.4}
	\left(\frac{2t}{e^{t}+1}\right)=\sum_{n=0}^{\infty}\frac{t^n}{n!}~\mathsf{G}_{n}, \,\,\,\,\,\,\, ( |t|<\pi; 1^{\alpha}:=1 ).
\end{equation}

The generating function of an Euler polynomials of order $\alpha$ (see\cite{Srivastva-garg}) is as folows:
\begin{equation}\label{e2.4}
	\sum_{n\geq0}\mathsf{E}_{n}^{(\alpha)}(x_{1})\frac{{t}^n}{n!}=e^{x_{1}{t}}\left(\frac{2}{e^{{t}}+1}\right)^{\alpha}\,\,\,\,\,\,\, ( |{t}|<\pi; 1^{\alpha}:=1 ).
\end{equation}

If we take $x_{1}=0$ in \eqref{e2.4} i.e. $\mathsf{E}_{n}^{(\alpha)}(0)=\mathsf{E}_{n}^{(\alpha)}$ are called Euler numbers defined as follows:
\begin{equation}\label{2.5}
	\sum_{n\geq0}\mathsf{E}_{n}^{(\alpha)}\frac{{t}^n}{n!}=\left(\frac{2}{e^{{t}}+1}\right)^{\alpha}.
\end{equation}

The generating function of an Bernaolli polynomials of order $\alpha$ (see \cite{Srivastva-garg})  is as follows:
\begin{equation}\label{2.44}
	\sum_{n\geq0}\mathsf{B}_{n}^{(\alpha)}(x_{1})\frac{{t}^n}{n!}=e^{x_{1}{t}}\left(\frac{t}{e^{{t}}-1}\right)^{\alpha}\,\,\,\,\,\,\, ( |{t}|<2\pi; 1^{\alpha}:=1 ).
\end{equation}

If we take $x_{1}=0$ in \eqref{2.44} i.e. $\mathsf{B}_{n}^{(\alpha)}(0)=\mathsf{B}_{n}^{(\alpha)}$ are called Bernaolli numbers defined as follows:
\begin{equation}\label{2.5}
	\sum_{n\geq0}\mathsf{B}_{n}^{(\alpha)}\frac{{t}^n}{n!}=\left(\frac{t}{e^{{t}}-1}\right)^{\alpha}.
\end{equation}

The generating function of an Genocchi polynomials of order $\alpha$ (see \cite{Srivastva-garg}) is as follows:
\begin{equation}\label{2.4}
\sum_{n\geq0}\mathsf{G}_{n}^{(\alpha)}(x_{1})\frac{{t}^n}{n!}=e^{x_{1}{t}}\left(\frac{2t}{e^{{t}}+1}\right)^{\alpha}\,\,\,\,\,\,\, ( |{t}|<\pi; 1^{\alpha}:=1 ).
\end{equation}

If we take $x_{1}=0$ in \eqref{2.4} i.e. $\mathsf{G}_{n}^{(\alpha)}(0)=\mathsf{G}_{n}^{(\alpha)}$ are called Genocchi numbers defined as follows:
\begin{equation}\label{2.5}
\sum_{n\geq0}\mathsf{G}_{n}^{(\alpha)}\frac{{t}^n}{n!}=\left(\frac{2t}{e^{{t}}+1}\right)^{\alpha}.
\end{equation}

\vspace{0.15cm} 
The generalized Apostol-Bernoulli polynomials $\mathsf{B}_{n}^{(\alpha)}(x_{1}, \lambda)$ of order $\alpha \in \mathbb{C}$ (see \cite{Luo-ming-sri}), defined by the following generating function to be

\begin{equation}\label{e1.17}
	e^{x_{1}{t}}\left(\frac{t}{\lambda e^{t}-1}\right)^{\alpha}=\sum_{n=0}^{\infty}\frac{t^n}{n!}~{\mathsf{B}}_{n}^{(\alpha)}(x_{1}; \lambda),\,\,\,\,(|t+ln{\lambda}|<2\pi; 1^{\alpha}:=1),
\end{equation}
with
$${\mathsf{B}}_{n}^{(\alpha)}(x_{1}; 1):={\mathsf{B}}_{n}^{(\alpha)}(x_{1})$$
and
$${\mathsf{B}}_{n}^{(\alpha)}(0; \lambda):={\mathsf{B}}_{n}^{(\alpha)}(\lambda)$$

which is known as Apostol-Bernoulli numbers ${\mathsf{B}}_{n}^{(\alpha)}(\lambda)$ of order $\alpha$.

\vspace{0.15cm}
The generalized Apostol-Euler polynomials $\mathsf{E}_{n}^{(\alpha)}(x_{1}, \lambda)$ of order $\alpha \in \mathbb{C}$ (see\cite{luo-ming}), defined by the following generating function to be

\begin{equation}\label{e1.18}
	e^{x_{1}{t}}\left(\frac{2}{\lambda e^{t}+1}\right)^{\alpha}=\sum_{n=0}^{\infty}\frac{t^n}{n!}~\mathsf{E}_{n}^{(\alpha)}(x_{1}; \lambda),\,\,\,\,(|t+ln{\lambda}|<\pi; 1^{\alpha}:=1),
\end{equation}
with
$$\mathsf{E}_{n}^{(\alpha)}(x_{1}; 1):=\mathsf{E}_{n}^{(\alpha)}(x_{1})$$
and
$$\mathsf{E}_{n}^{(\alpha)}(0; \lambda):=\mathsf{E}_{n}^{(\alpha)}(\lambda)$$

which is known as Apostol-Euler numbers of order $\alpha$.

\vspace{0.15cm}
The generalized Apostol-Genocchi polynomials $\mathsf{G}_{n}^{(\alpha)}(x_{1}, \lambda)$  of order $\alpha \in \mathbb{C}$ (see \cite{luo-ming2}), defined by the following generating function to be

\begin{equation}\label{e1.19}
	e^{x_{1}{t}}\left(\frac{2t}{\lambda e^{t}+1}\right)^{\alpha}=\sum_{n=0}^{\infty}\frac{t^n}{n!}~\mathsf{G}_{n}^{(\alpha)}(x_{1}; \lambda),\,\,\,\,(|t+ln{\lambda}|<\pi; 1^{\alpha}:=1),
\end{equation}
with
$$\mathsf{G}_{n}^{(\alpha)}(x_{1}; 1):=\mathsf{G}_{n}^{(\alpha)}(x_{1})$$
and
$$\mathsf{G}_{n}^{(\alpha)}(0; \lambda):=\mathsf{G}_{n}^{(\alpha)}(\lambda)$$

The second kind of Stirling polynomials and Stirling numbers are denoted by  $\mathcal{S}_2(m,n;x_{1})$, $\mathcal{S}_2(m,n)$ and defined by following generating functions (see\cite{Bell,Boas}):
\begin{equation}\label{2.6}
	\sum_{n\geq0}\mathcal{S}_{2}(m,n;x_{1})\frac{{t}^m}{m!}=\frac{(e^{{t}}-1)^{n}}{n!}e^{{t}x_{1}},
\end{equation}

when $x_{1}=0$ in \eqref{2.6} i.e. $\mathcal{S}_2(m,n;0)=\mathcal{S}_2(m,n)$ are called Stirling number and  generating functions define as follows (see\cite{Bell,Boas}):
\begin{equation}\label{2.7}
	\sum_{n\geq0}\mathcal{S}_{2}(m,n)\frac{{t}^m}{m!}=\frac{(e^{{t}}-1)^{n}}{n!}.
\end{equation}

In 2013, Lu and Luo (\cite{lu-luo}) introduced the generalized Apostol type polynomials of order $\alpha$ defined by the following generating function:

\begin{equation}\label{e1.20}
	e^{x_{1}{t}}\left(\frac{2^{\eta}t^{\delta}}{\lambda e^{t}+1}\right)^{\alpha}=\sum_{n=0}^{\infty}\frac{t^n}{n!}~\mathcal{F}_{n}^{(\alpha)}(x_{1}; \lambda;\eta,\delta),\,\,\,\,(|t|<|log(-\lambda)|),
\end{equation}
$$(\alpha\in \mathbb{N}_{0}, \lambda, \eta, \delta \in \mathbb{C}).$$
Where 
\begin{equation}\label{1.21}
\mathcal{F}_{n}^{(\alpha)}(\lambda;\eta,\delta)=\mathcal{F}_{n}^{(\alpha)}(0; \lambda;\eta,\delta,
\end{equation}
are called the Apostol numbers of order $\alpha$. Noted that equation \eqref{e1.20} can be reduced to equation \eqref{e1.17}, \eqref{e1.18} and \eqref{e1.19} as follows:

\begin{equation}\label{e1.22}
	\mathsf{B}_{n}^{(\alpha)}(x_{1};\lambda)= (-1)^{\alpha}\mathcal{F}(x_{1};-\lambda;0,1).
\end{equation}

\begin{equation}\label{e1.23}
	\mathsf{E}_{n}^{(\alpha)}(x_{1};\lambda)= \mathcal{F}(x_{1};\lambda;1,0).
\end{equation}

\begin{equation}\label{e1.24}
	\mathsf{G}_{n}^{(\alpha)}(x_{1};\lambda)= \mathcal{F}(x_{1};\lambda;1,1).
\end{equation}

Inspired by the result obtained in (see \cite{Araci-Acikgoz,khan-husain,lu-luo}). In the present article, we mentioned a new type of mix type polynomials known as Bell based Apostol polynomials of order $\alpha$ and studied their some special cases and  elementary properties.

\section{\bf Bell based Apostol type polynomials ${_{\mathfrak{B}} \mathcal{F} _{n}^{(\alpha)}}(x_{1},x_{2};\lambda;\eta,\delta)$}
In this section, we introduce the Bell based Apostol type polynomials of order $\alpha$ and investigate their numerous relation such as correlation formulae, implicit summation formulae, derivative formulae.

\begin{definition}
The generating function for Bell based Apostol type polynomials of order $\alpha$ define as follows:

For any $n\in\mathbb{N}$, $\alpha \in \mathbb{N}_{0}$ and $\lambda, \eta, \delta \in \mathbb{C}$, we define Bell based Apostol type polynomials of order $\alpha$ as:\\
\begin{equation}\label{3.1}
\sum_{n\geq0}{_{\mathfrak{B}} \mathcal{F} _{n}^{(\alpha)}}(x_{1},x_{2};\lambda;\eta,\delta)\frac{{t}^n}{n!}=\left(\frac{2^{\eta}t^{\delta}}{\lambda e^{{t}}+1}\right)^{\alpha}e^{{x_{1}{t}}+x_{2}(e^t -1)},\,\,\,\,\,\,\, (|t|<|log(-\lambda)|).
\end{equation}
\end{definition}

If $x_{1}=0$ and $x_{2}=1$ in \eqref{3.1} then we get a Bell Based Apostol type numbers of order $\alpha$, which is defined  as follows:
\begin{equation}\label{3.2}
	\sum_{n\geq0}{_{\mathfrak{B}} \mathcal{F} _{n}^{(\alpha)}}(0,1;\lambda;\eta,\delta)\frac{{t}^n}{n!}=\left(\frac{2^{\eta}t^{\delta}}{\lambda e^{{t}}+1}\right)^{\alpha}e^{(e^t -1)},\,\,\,\,\,\,\, (|t|<|log(-\lambda)|).
\end{equation}

\begin{remark}
	If we choose $\eta=\delta=\lambda=1$ in \eqref{3.1}, then we have to reduces Bell based Apostol polynomilas of order $\alpha$ into Bell based Euler polynomials of order $\alpha$ (see \cite{khan-husain}).
\end{remark}

\begin{remark}
	If we choose $\alpha=0$ in \eqref{3.1}, we have to reduce Bell based Apostol type polynomials of order $\alpha$, into bivariate Bell polynomials (see\cite{Araci-Acikgoz}) defined as follows:
	\begin{equation*}
	\sum_{n\geq0}{_{\mathfrak{B}} \mathcal{F} _{n}^{(0)}}(x_{1},x_{2};\lambda;\eta,\delta)\frac{{t}^n}{n!}=e^{x_{1}t+x_{2}(e^{t} -1)}.
	\end{equation*}
\end{remark}

\begin{remark}
In case $x_{2}=0$ in \eqref{3.1} the Bell based Apostol type polynomials of order $\alpha$ reduces to the Apostol type polynomials (see\cite{lu-luo}) of order $\alpha$ defined as follows:
\begin{equation*}
\sum_{n\geq0}{_{\mathfrak{B}} \mathcal{F} _{n}^{(\alpha)}}(x_{1},0;\lambda;\eta,\delta)\frac{{t}^n}{n!}=\left(\frac{2^{\eta}t^{\delta}}{\lambda e^{{t}}+1}\right)^{\alpha}e^{xt}.
\end{equation*}
\end{remark}

\subsection{Some elementary properties of ${_{\mathfrak{B}} \mathcal{F} _{n}^{(\alpha)}}(x_{1},x_{2};\lambda;\eta,\delta)$ }

 \begin{theorem}
 For $\alpha \in \mathbb{N}_{0}$, $\lambda, \eta, \delta \in \mathbb{C}$ and $n \in \mathbb{N}$ the following relation of Bell based Apostol type polynomials hold true:
 
 \begin{equation}\label{3.3}
 	{_{\mathfrak{B}} \mathcal{F} _{n}^{(\alpha)}}(x_{1},x_{2};\lambda;\eta,\delta)=\sum\limits_{k=0}^{n}\binom{n}{k} {_{\mathfrak{B}} \mathcal{F} _{k}^{(\alpha)}}(x_{1};\lambda;\eta,\delta) {\mathfrak{B}}_{n-k}(x_{2}).
 \end{equation}

\begin{proof} By using the relation \eqref{3.1}, we have
	
\begin{equation*}
\aligned
\sum_{n\geq0}{_{\mathfrak{B}} \mathcal{F} _{n}^{(\alpha)}}(x_{1},x_{2};\lambda;\eta,\delta)\frac{{t}^n}{n!}=&\left(\frac{2^{\eta}t^{\delta}}{\lambda e^{{t}}+1}\right)^{\alpha}e^{{x_{1}{t}}+x_{2}(e^{t} -1)}\\
=&\left\{ \left(\frac{2^{\eta}t^{\delta}}{\lambda e^{{t}}+1}\right)^{\alpha}e^{x_{1}{t}}\right\} \left\{e^{x_{2}(e^{t} -1)}\right\}\\
=& \left\{\sum_{k\geq0}{\mathcal{F} _{k}^{(\alpha)}}(x_{1};\lambda;\eta,\delta)\frac{{t}^k}{k!}\right\}\left\{\sum_{n\geq0}{\mathfrak{B}_n}(x_{2})\frac{{t}^n}{n!}\right\}.
\endaligned
\end{equation*}

Now, using the series rearrangement technique, we get
\begin{equation*}
\aligned
\sum_{n\geq0}{_{\mathfrak{B}} \mathcal{F} _{n}^{(\alpha)}}(x_{1},x_{2};\lambda;\eta,\delta)\frac{{t}^n}{n!}=& \sum_{n\geq0}\left\{\sum\limits_{k=0}^{n}\binom{n}{k} { \mathcal{F} _{k}^{(\alpha)}}(x_{1};\lambda;\eta,\delta) {\mathfrak{B}}_{n-k}(x_{2})\right\}\frac{{t}^n}{n!}.
\endaligned
\end{equation*}

By equating the same power of t both side, we obtained the result \eqref{3.3}. 
\end{proof}
\end{theorem}

\begin{theorem}
For any $\alpha \in \mathbb{N}_{0}$, $\lambda, \eta, \delta \in \mathbb{C}$ and $n \in \mathbb{N}$ the following relation hold true:
\begin{equation}\label{3.4}
	{_{\mathfrak{B}} \mathcal{F} _{n}^{(\alpha)}}(x_{1},x_{2};\lambda;\eta,\delta)=\sum\limits_{k=0}^{n}\binom{n}{k} { \mathcal{F} _{k}^{(\alpha)}}(0;\lambda;\eta,\delta) {\mathfrak{B}}_{n-k}(x_{1};x_{2}).
\end{equation}

\begin{proof}
 By	using the generating function \eqref{3.1}, we have
 \begin{equation*}
 	\aligned
 	\sum_{n\geq0}{_{\mathfrak{B}} \mathcal{F} _{n}^{(\alpha)}}(x_{1},x_{2};\lambda;\eta,\delta)\frac{{t}^n}{n!}=&\left(\frac{2^{\eta}t^{\delta}}{\lambda e^{{t}}+1}\right)^{\alpha}e^{{x_{1}{t}}+x_{2}(e^{t} -1)}\\
 	=&\left\{ \left(\frac{2^{\eta}t^{\delta}}{\lambda e^{{t}}+1}\right)^{\alpha}\right\} \left\{e^{x_{1}{t}+x_{2}(e^{t} -1)}\right\}\\
 	=& \left\{\sum_{k\geq0}{ \mathcal{F} _{k}^{(\alpha)}}(0;\lambda;\eta,\delta)\frac{{t}^k}{k!}\right\}\left\{\sum_{n\geq0}{\mathfrak{B}_n}(x_{1};x_{2})\frac{{t}^n}{n!}\right\}.\\
 	\endaligned
 \end{equation*}

After applying the series rearrangement technique, we obtained the result \eqref{3.4}.
\end{proof}
\end{theorem}	
	
\begin{theorem}
For $\alpha \in \mathbb{N}_{0}$, $\lambda, \eta, \delta \in \mathbb{C}$ and $n \in \mathbb{N}$, then the following relation of the Bell based Apostol type polynomials hold true:	
\begin{equation}\label{3.5}
	{_{\mathfrak{B}} \mathcal{F} _{n}^{(\alpha)}}(x_{1},x_{2};\lambda;\eta,\delta)=\sum\limits_{k=0}^{n}\binom{n}{k} {_{\mathfrak{B}} \mathcal{F} _{k}^{(\alpha)}}(x_{2};\lambda;\eta,\delta)\,\,{x_{1}}^{n-k}.
\end{equation}
\begin{proof}
 Using the relation \eqref{3.1}, we have
 \begin{equation*}
 	\aligned
 {_{\mathfrak{B}} \mathcal{F} _{n}^{(\alpha)}}(x_{1},x_{2};\lambda;\eta,\delta)\frac{{t}^n}{n!}=&\left(\frac{2^{\eta}t^{\delta}}{\lambda e^{{t}}+1}\right)^{\alpha}e^{{x_{1}{t}}+x_{2}(e^{t} -1)}\\
 	=&\left\{ \left(\frac{2^{\eta}t^{\delta}}{\lambda e^{{t}}+1}\right)^{\alpha}e^{x_{2}(e^{{t}}-1)}\right\} \left\{e^{x_{1}{t}}\right\}\\
 	=& \left\{\sum_{k\geq0} {_{\mathfrak{B}} \mathcal{F} _{k}^{(\alpha)}}(x_{2};\lambda;\eta,\delta)\frac{{t}^k}{k!}\right\}\left\{\sum_{n\geq0}\frac{(x_{1}{t})^n}{n!}\right\}\\
 	=& \left\{\sum_{n\geq0} \sum_{k\geq0} {_{\mathfrak{B}} \mathcal{F} _{k}^{(\alpha)}}(x_{2};\lambda;\eta,\delta)\,\,\frac{{x_{1}}^n}{n!}\frac{{t}^{n+k}}{k!}\right\},
 	\endaligned
 \end{equation*}

using series rearrangement technique, we get 
\begin{equation*}
 	\sum_{n\geq0}{_{\mathfrak{B}} \mathcal{F} _{n}^{(\alpha)}}(x_{1},x_{2};\lambda;\eta,\delta)\frac{{t}^n}{n!}= \sum_{n\geq0}\left\{ \sum_{k=0}^{n}\binom{n}{k}{_{\mathfrak{B}} \mathcal{F} _{k}^{(\alpha)}}(x_{2};\lambda;\eta,\delta)\,\,{x_{1}}^{n-k}\right\}\frac{{t}^{n}}{n!}.
\end{equation*}

By equating the same power of t both side, we obtained the result \eqref{3.5}. 
\end{proof}
\end{theorem}

\begin{theorem}
For $\alpha \in \mathbb{N}_{0}$, $\lambda, \eta, \delta \in \mathbb{C}$ and $n \in \mathbb{N}$, then the following relation of the Bell based Apostol type polynomials hold true:
	\begin{equation}\label{3.6}
		{_{\mathfrak{B}} \mathcal{F} _{n}^{(\alpha)}}(x_{1}+x_{2},x_{3};\lambda;\eta,\delta)=\sum\limits_{k=0}^{n}\binom{n}{k} { \mathcal{F} _{k}^{(\alpha)}}(x_{1};\lambda;\eta,\delta) {\mathfrak{B}}_{n-k}(x_{2},x_{3}).
	\end{equation}
	
	\begin{proof}
		By	using the generating function \eqref{3.1}, we have
		\begin{equation*}
			\aligned
			\sum_{n\geq0}{_{\mathfrak{B}} \mathcal{F} _{n}^{(\alpha)}}(x_{1}+x_{2},x_{3};\lambda;\eta,\delta)\frac{{t}^n}{n!}=&\left(\frac{2^{\eta}t^{\delta}}{\lambda e^{{t}}+1}\right)^{\alpha}e^{(x_{1}+x_{2}){t}+x_{3}(e^{t} -1)}\\
			=&\left\{ \left(\frac{2^{\eta}t^{\delta}}{\lambda e^{{t}}+1}\right)^{\alpha}e^{x_{1}t}\right\} \left\{e^{x_{2}{t}+x_{3}(e^{t} -1)}\right\}\\
			=& \left\{\sum_{k\geq0}{ \mathcal{F} _{k}^{(\alpha)}}(x_{1};\lambda;\eta,\delta)\frac{{t}^k}{k!}\right\}\left\{\sum_{n\geq0}{\mathfrak{B}_n}(x_{2},x_{3})\frac{{t}^n}{n!}\right\}.\\
			\endaligned
		\end{equation*}
		
		After applying series rearrangement technique, we obtained the result \eqref{3.6}.
	\end{proof}
\end{theorem}

\section{\bf Implicit summation formulas of ${_{\mathfrak{B}} \mathcal{F} _{n}^{(\alpha)}}(x_{1},x_{2};\lambda;\eta,\delta)$}
In this section, we have to discuss useful identities such as the implicit summation formula for the Bell Based Apostol type polynomials of order $\alpha$, which is define in following theorems as follows:

\begin{theorem}
For $\alpha_{1}, \alpha_{2} \in \mathbb{N}_{0}$, $\lambda, \eta, \delta \in \mathbb{C}$ and $n \in \mathbb{N}$, then the following relation of Bell based Apostol type polynomials hold true::

$${_{\mathfrak{B}} \mathcal{F}_{n}^{(\alpha_{1}+\alpha_{2})}}(x_{1}+x_{2},y_{1}+y_{2};\lambda;\eta,\delta)$$	
 \begin{equation}\label{4.1}
\aligned
=&\sum\limits_{k=0}^{n}\binom{n}{k} {_\mathfrak{B} {\mathcal{F}}_{k}^{(\alpha_{1})}}(x_{1},y_{1};\lambda;\eta,\delta)\,\, {_\mathfrak{B}}{\mathcal{F}}^{(\alpha_{2})}_{n-k}(x_{2},y_{2};\lambda;\eta,\delta).
\endaligned
 \end{equation}	
	
\begin{proof}
By using the following identity

$$\left(\frac{2^{\eta}t^{\delta}}{\lambda e^{{t}}+1}\right)^{\alpha_{1}+\alpha_{2}}e^{{(x_{1}+x_{2})t}+(y_{1}+y_{2})(e^t -1)}$$

$$=\left\{ \left(\frac{2^{\eta}t^{\delta}}{\lambda e^{{t}}+1}\right)^{\alpha_{1}}e^{x_{1}{t}+y_{1}(e^{t} -1)}\right\} \left\{ \left(\frac{2^{\eta}t^{\delta}}{\lambda e^{{t}}+1}\right)^{\alpha_{2}}e^{x_{2}{t}+y_{2}(e^{t} -1)}\right\}.$$

 Applying the above identity in the generating function \eqref{3.1}, we have
	
$$\sum_{n\geq0}{_{\mathfrak{B}} \mathcal{F}_{n}^{(\alpha_{1}+\alpha_{2})}}(x_{1}+x_{2},y_{1}+y_{2};\lambda;\eta;\delta)\frac{{t}^n}{n!}$$
\begin{equation*}
\aligned	
=&\left(\frac{2^{\eta}t^{\delta}}{\lambda e^{{t}}+1}\right)^{\alpha_{1}+\alpha_{2}}e^{{(x_{1}+x_{2}){t}}+(y_{1}+y_{2})(e^{t} -1)}\\
=&\left\{ \left(\frac{2^{\eta}t^{\delta}}{\lambda e^{{t}}+1}\right)^{\alpha_{1}}e^{x_{1}{t}+y_{1}(e^{t} -1)}\right\} \left\{ \left(\frac{2^{\eta}t^{\delta}}{\lambda e^{{t}}+1}\right)^{\alpha_{2}}e^{x_{2}{t}+y_{2}(e^{t} -1)}\right\}\\
=&\left\{\sum_{k\geq0}{_\mathfrak{B}}{\mathcal{F}_{k}^{(\alpha_{1})}}(x_{1},y_{1};\lambda;\eta,\delta)\frac{{t}^k}{k!}\right\} \left\{\sum_{n\geq0}{_\mathfrak{B}}{\mathcal{F}_{n}^{(\alpha_{2})}}(x_{2},y_{2};\lambda;\eta,\delta)\frac{{t}^n}{n!}\right\}\\
=&\left\{\sum_{n\geq0}\sum_{k\geq0}{_\mathfrak{B}}{\mathcal{F}_{k}^{(\alpha_{1})}}(x_{1},y_{1};\lambda;\eta,\delta)\,\,{_\mathfrak{B}}{\mathcal{F}_{n}^{(\alpha_{2})}}(x_{2},y_{2};\lambda;\eta,\delta)\,\,\frac{{t}^{n+k}}{n! k!}\right\},
\endaligned	
\end{equation*}

using the series rearrangement technique, we obtain
$$\sum_{n\geq0}{_{\mathfrak{B}} \mathcal{F}_{n}^{(\alpha_{1}+\alpha_{2})}}(x_{1}+x_{2},y_{1}+y_{2};\lambda;\eta,\delta)\frac{{t}^n}{n!}$$
\begin{equation*}
\aligned
=&\sum_{n\geq0}\left\{\sum_{k=0}^{n}\binom{n}{k}{_\mathfrak{B}}{\mathcal{F}_{k}^{(\alpha_{1})}}(x_{1},y_{1};\lambda;\eta,\delta)\,\,{_\mathfrak{B}}{\mathcal{F}_{n-k}^{(\alpha_{2})}}(x_{2},y_{2};\lambda;\eta,\delta)\,\,\right\}\frac{{t}^{n}}{n!}.
\endaligned
\end{equation*}

Now, equating the same power of ${t}$ both side, we obtained the result \eqref{4.1}.
\end{proof}	
\end{theorem}

\begin{remark}
In case, if we take $\alpha_{1}=\alpha$, $\alpha_{2}=0$, $x_{1}=x$, $x_{2}=1$, $y_{1}=y$ and $y_{2}=0$ in \eqref{4.1}, we have
\begin{equation}\label{4.2}
	{_{\mathfrak{B}} \mathcal{F}_{n}^{(\alpha)}}(x+1,y;\lambda;\eta,\delta)=\sum\limits_{k=0}^{n}\binom{n}{k} {_\mathfrak{B} {\mathcal{F}}_{k}^{(\alpha)}}(x,y;\lambda;\eta,\delta),
\end{equation}
 \end{remark}

\begin{theorem}
For $\alpha \in \mathbb{N}_{0}$, $\lambda, \eta, \delta \in \mathbb{C}$ and $n \in \mathbb{N}$, then the following relation of the Bell based Apostol type polynomials hold true:

\begin{equation}\label{4.4}
{_{\mathfrak{B}} \mathcal{F}_{n+1}^{(\alpha)}}(x_{1}+1,x_{2};\lambda;\eta,\delta)-	{_{\mathfrak{B}} \mathcal{F}_{n+1}^{(\alpha)}}(x_{1},x_{2};\lambda;\eta,\delta)=\sum\limits_{k=0}^{n}\binom{n+1}{k} {_\mathfrak{B} {\mathcal{F}}_{k}^{(\alpha)}}(x_{1},x_{2};\lambda;\eta,\delta).
\end{equation}

\begin{proof}
Using the relation \eqref{3.1}, we get
$$\sum_{n\geq0}{_{\mathfrak{B}} \mathcal{F}_{n}^{(\alpha)}}(x_{1}+1,x_{2};\lambda;\eta,\delta)\frac{{t}^n}{n!}-\sum_{n\geq0}{_{\mathfrak{B}} \mathcal{F}_{n}^{(\alpha)}}(x_{1},x_{2};\lambda;\eta,\delta)\frac{{t}^n}{n!}$$
\begin{equation*}
\aligned
=&\left(\frac{2^{\eta}t^{\delta}}{`\lambda e^{{t}}+1}\right)^{\alpha}e^{{(x_{1}+1){t}}+x_{2}(e^{t} -1)}-\left(\frac{2^{\eta}t^{\delta}}{\lambda e^{{t}}+1}\right)^{\alpha}e^{{x_{1}{t}}+x_{2}(e^{t} -1)}\\
=& \left(\frac{2^{\eta}t^{\delta}}{\lambda e^{{t}}+1}\right)^{\alpha}e^{{x_{1}{t}}+x_{2}(e^{t} -1)}\,\,(e^{{t}}-1)\\
=& \left\{\sum_{k\geq0}{_{\mathfrak{B}} \mathcal{F}_{k}^{(\alpha)}}(x_{1},x_{2};\lambda;\eta,\delta)\frac{{t}^k}{k!}\right\}\,\,\left\{\sum_{n\geq0}\frac{{t}^{n+1}}{(n+1)!}\right\}.
\endaligned
\end{equation*}

After, applying the series rearrangement technique, we obtain the result \eqref{4.4} 
\end{proof}
\end{theorem}

\begin{theorem}
For $\alpha \in \mathbb{N}_{0}$, $\lambda, \eta, \delta \in \mathbb{C}$ and $n \in \mathbb{N}$, then the following relation of Bell based Apostol type polynomials hold true:
\begin{equation}\label{4.5}
	{\mathfrak{B}_n}(x_{1};x_{2})=\frac{n!\left(\lambda\,{_{\mathfrak{B}} \mathcal{F}_{n+\delta}^{(1)}}(x_{1}+1,x_{2};\lambda;\eta,\delta)+	{_{\mathfrak{B}} \mathcal{F}_{n+\delta}^{(1)}}(x_{1},x_{2};\lambda;\eta,\delta)\right)}{2^{\eta}(n+\delta)!}.
\end{equation}
\begin{proof}
 Using the generating function \eqref{3.1} for $\alpha=1$ and definition of the bivariate Bell polynomials, we get
 $$\sum_{n\geq0}{\mathfrak{B}_n}(x_{1};x_{2})\frac{{t}^n}{n!}$$
 \begin{equation*}
 	\aligned
 	=& e^{x_{1}{t}+{x_{2}(e^{{t}}-1)}}\\
 	=&\frac{\lambda e^{{t}}+1}{2^{\eta}t^{\delta}}\,\,\left\{\sum_{k=0}^{\infty}{_{\mathfrak{B}} \mathcal{F}_{k}^{(1)}}(x_{1},x_{2};\lambda;\eta,\delta)\right\}\\
 	=& \frac{\lambda e^{{t}}+1}{2^{\eta}t^{\delta}}\left\{\left(\frac{2^{\eta}t^{\delta}}{\lambda e^{{t}}+1}\right)e^{{x_{1}{t}}+x_{2}(e^{t} -1)}\right\}\\
 	=&\frac{1}{2^{\eta}t^{\delta}}\left\{\lambda\left(\frac{2^{\eta}t^{\delta}}{\lambda e^{{t}}+1}\right)e^{{(x_{1}+1){t}}+x_{2}(e^{t} -1)}+\left(\frac{2^{\eta}t^{\delta}}{\lambda e^{{t}}+1}\right)e^{{x_{1}{t}}+x_{2}(e^{t} -1)}\right\}\\
 	=&\frac{1}{2^{\eta}t^{\delta}}\left\{\lambda\sum_{n\geq0}{_{\mathfrak{B}}\mathcal{F}_{n}^{(1)}}(x_{1}+1,x_{2};\lambda;\eta,\delta)\frac{{t}^n}{n!}+\sum_{n\geq0}{_{\mathfrak{B}} \mathcal{F}_{n}^{(1)}}(x_{1},x_{2};\lambda;\eta,\delta)\frac{{t}^n}{n!}\right\}\\
 	=&\frac{1}{2^{\eta}}\left\{\lambda\sum_{n\geq0}{_{\mathfrak{B}}\mathcal{F}_{n}^{(1)}}(x_{1}+1,x_{2};\lambda;\eta,\delta)\frac{{t}^{n-\delta}}{n!}+\sum_{n\geq0}{_{\mathfrak{B}} \mathcal{F}_{n}^{(1)}}(x_{1},x_{2};\lambda;\eta,\delta)\frac{{t}^{n-\delta}}{n!}\right\}\\
 	=&\frac{n!}{2^{\eta}(n+\delta)!}\left\{\lambda\sum_{n\geq0}{_{\mathfrak{B}}\mathcal{F}_{n+\delta}^{(1)}}(x_{1}+1,x_{2};\lambda;\eta,\delta)\frac{{t}^n}{n!}+\sum_{n\geq0}{_{\mathfrak{B}} \mathcal{F}_{n+\delta}^{(1)}}(x_{1},x_{2};\lambda;\eta,\delta)\frac{{t}^n}{n!}\right\}.
 	\endaligned
 \end{equation*}
By equating the same power of ${t}$ both side, we obtained the result \eqref{4.5}.
\end{proof}
\end{theorem}

\begin{theorem}
For $\alpha \in \mathbb{N}_{0}$, $\lambda, \eta, \delta \in \mathbb{C}$ and $n \in \mathbb{N}$, then the following relation of Bell based Apostol type polynomials hold true:
\begin{equation}\label{4.7}
{_{\mathfrak{B}} \mathcal{F}_{n}^{(\alpha)}}(x_{1},x_{2};\lambda;\eta,\delta)= \sum_{j=0}^{n}\sum_{k\geq0}\binom{n}{j}(x_{1})_{k}\,\mathcal{S}_{2}(j, k)\,{_{\mathfrak{B}} \mathcal{F}_{n}^{(\alpha)}}(x_{2};\lambda;\eta,\delta).
\end{equation}

\begin{proof}
 By using the relation \eqref{3.1}, we have
 \begin{equation*}
 	\aligned
 	\sum_{n\geq0}{_{\mathfrak{B}} \mathcal{F}_{n}^{(\alpha)}}(x_{1},x_{2};\lambda;\eta,\delta)\frac{{t}^n}{n!}=&\left(\frac{2^{\eta}t^{\delta}}{\lambda e^{{t}}+1}\right)^{\alpha}e^{{x_{1}{t}}+x_{2}(e^{t} -1)}\\
 	=&\left(\frac{2^{\eta}t^{\delta}}{\lambda e^{{t}}+1}\right)^{\alpha}e^{x_{2}(e^{t} -1)}\,e^{x_{1}{t}}\\
 	=&\left(\frac{2^{\eta}t^{\delta}}{\lambda e^{{t}}+1}\right)^{\alpha}e^{x_{2}(e^{t} -1)}\,(1+e^{{t}} -1)^{x_{1}}\\
 	=&\left\{\sum_{n\geq0}{_{\mathfrak{B}} \mathcal{F}_{n}^{(\alpha)}}(x_{2};\lambda;\eta,\delta)\frac{{t}^n}{n!}\right\}\left\{\sum_{k\geq0}(x_{1})_{k}\,\frac{(e^{{t}} -1)^{k}}{k!}\right\}\\
 	=&\left\{\sum_{n\geq0}{_{\mathfrak{B}} \mathcal{F}_{n}^{(\alpha)}}(x_{2};\lambda;\eta,\delta)\frac{{t}^n}{n!}\right\}\left\{\sum_{k\geq0}(x_{1})_{k}\,\sum_{j\geq0}\mathcal{S}_{2}(j, k)\frac{{t}^{j}}{j!}\right\}.
 	\endaligned
 \end{equation*}
By using the series rearrangement technique, we obtained the result \eqref{4.7}.
\end{proof}
\end{theorem}

\section{\bf Derivative formulas of ${_{\mathfrak{B}} \mathcal{F} _{n}^{(\alpha)}}(x_{1},x_{2};\lambda;\eta,\delta)$}
\begin{theorem}
For $\alpha \in \mathbb{N}_{0}$, $\lambda, \eta, \delta \in \mathbb{C}$ and $n \in \mathbb{N}$, then the following difference operator formula of Bell based Apostol type polynomials w.r.t. $x_{1}$ hold true:
\begin{equation}\label{5.1}
	\frac{\partial}{\partial x_{1}}\,\, {_{\mathfrak{B}} \mathcal{F}_{n}^{(\alpha)}}(x_{1},x_{2};\lambda;\eta,\delta)= n\,\,{_{\mathfrak{B}} \mathcal{F}_{n-1}^{(\alpha)}}(x_{1},x_{2};\lambda;\eta,\delta),
\end{equation}
which hold for all $n\in\mathbb{N}$.

\begin{proof}
  we know that
  \begin{equation}\label{5.2}
  \frac{\partial}{\partial x_{1}}\,\,{e^{x_{1}{t}+x_{2}(e^{{t}}-1)}}=t\,\,e^{x_{1}{t}+x_{2}(e^{{t}}-1)}.
  \end{equation}

  By using the definition \eqref{3.1} in \eqref{5.2}, we obtained the result \eqref{5.1}.
\end{proof}
\end{theorem}

\begin{theorem}
For $\alpha \in \mathbb{N}_{0}$, $\lambda, \eta, \delta \in \mathbb{C}$ and $n \in \mathbb{N}$, then the following difference operator formula of Bell based Apostol type polynomials w.r.t. $x_{2}$ hold true:
\begin{equation}\label{5.3}
	\frac{\partial}{\partial x_{2}}\,\, {_{\mathfrak{B}} \mathcal{F}_{n}^{(\alpha)}}(x_{1},x_{2};\lambda;\eta,\delta)=\left\{{_{\mathfrak{B}} \mathcal{F}_{n}^{(\alpha)}}(x_{1}+1,x_{2};\lambda;\eta,\delta)-{_{\mathfrak{B}} \mathcal{F}_{n}^{(\alpha)}}(x_{1},x_{2};\lambda;\eta,\delta)\right\},
\end{equation}
which hold for all $n\in\mathbb{N}$.
\begin{proof}
	By using well known derivative properties
	\begin{equation}\label{5.4}
		\frac{\partial}{\partial x_{2}}\,\,{e^{x_{1}{t}+x_{2}(e^{{t}}-1)}}=(e^{{t}}-1)\,\,e^{x_{1}{t}+x_{2}(e^{{t}}-1)}.
	\end{equation}

Now, using the definition \eqref{3.1}, we get
\begin{equation*}
\frac{\partial}{\partial x_{2}}\left\{\sum_{n\geq0}{_{\mathfrak{B}} \mathcal{F}_{n}^{(\alpha)}}(x_{1},x_{2};\lambda;\eta,\delta)\frac{{t}^n}{n!}\right\}
\end{equation*}

 \begin{equation*}
\aligned
=&\frac{\partial}{\partial x_{2}}\left\{\left(\frac{2^{\eta}t^{\delta}}{\lambda e^{{t}}+1}\right)^{\alpha}e^{{x_{1}{t}}+x_{2}(e^{t} -1)}\right\}\\
=&\left\{ \left(\frac{2^{\eta}t^{\delta}}{\lambda e^{{t}}+1}\right)^{\alpha}e^{x_{1}{t}+x_{2}(e^{{t}}-1)}\right\}(e^{{t}}-1)\\
=&\left\{\left(\frac{2^{\eta}t^{\delta}}{\lambda e^{{t}}+1}\right)^{\alpha}e^{(x_{1}+1){t}+x_{2}(e^{{t}}-1)}\right\}-\left(\frac{2^{\eta}t^{\delta}}{\lambda e^{{t}}+1}\right)^{\alpha}e^{x_{1}{t}+x_{2}(e^{{t}}-1)}\\
=&\left\{\sum_{n\geq0}{_{\mathfrak{B}}\mathcal{F}_{n}^{(\alpha)}}(x_{1}+1,x_{2};\lambda;\eta,\delta)\frac{{t}^n}{n!}\right\}-\left\{\sum_{n\geq0}{_{\mathfrak{B}}\mathcal{F}_{n}^{(\alpha)}}(x_{1},x_{2};\lambda;\eta,\delta)\frac{{t}^{n}}{n!}\right\}\\
=&\left\{\sum_{n\geq0}{_{\mathfrak{B}}\mathcal{F}_{n}^{(\alpha)}}(x_{1}+1,x_{2};\lambda;\eta,\delta)-\sum_{n\geq0}{_{\mathfrak{B}}\mathcal{F}_{n}^{(\alpha)}}(x_{1},x_{2};\lambda;\eta,\delta)\right\}\frac{{t}^n}{n!}.
\endaligned
\end{equation*}

By equating  the same power of ${t}$ both side, we obtained the result \eqref{5.3}.
\end{proof}	
\end{theorem}

\section{\bf Conclusions}
Motivated by numerous applications in various field of Mathematical sciences such as  Number theory and Combinatorial analysis etc. In this paper, we introduced a mix type Bell based Apostol type polynomials of order $\alpha$ and studies their various elementary properties and identities such as correlation formulas, Implicit summation formulas and derivative formulas. In result obtained in this paper very general and specialized to yield a large number of new and known identities involving basic and unified polynomials given by other authors. 


\end{document}